\documentclass[12pt]{article}

\usepackage{amsmath, amssymb}
\usepackage[all]{xy}
\usepackage{color}
\usepackage{graphicx}

\newtheorem{theorem}{Theorem}[section]
\newtheorem{lemma}[theorem]{Lemma}

\newtheorem{corollary}[theorem]{Corollary}

\newcommand{\diff}{\text{\sf Diff}}

\newcommand{\n}{\noindent}
\newcommand{\so}{{}_{SO(3)}}

\newcommand{\F}{\mathbf F}
\newcommand{\R}{\mathbb R}
\newcommand{\Z}{\mathbb Z}
\newcommand{\K}{\mathbb K}
\newcommand{\mcg}{\Gamma}
\newcommand{\rp}{\mathbb R {\rm P}}

\newcommand{\Conf}[1]{F_k(\R{\rm P}^{#1})}
\newcommand{\conf}[1]{F_k(\R{\rm P}^{#1})   /  \Sigma_k}

\begin{document}
\title{On the cohomology of the mapping class group \\ of the punctured projective plane}

\author{Miguel A. Maldonado \\ {\small Unidad Acad\'emica de Matem\'aticas, Universidad Aut\'onoma}\\ 
{\small de Zacatecas, Zacatecas 98000, MEXICO}\\
{\small mmaldonado@matematicas.reduaz.mx}
\and 
Miguel A. Xicot\'encatl \\ {\small Departamento de Matem\'aticas, Centro de Investigaci\'on y de} \\
{\small Estudios Avanzados del IPN. Mexico City 07360, MEXICO}\\
{\small xico@math.cinvestav.mx}
}

\date{}

\maketitle

\begin{abstract}
The mapping class group $\Gamma^k(N_g)$ of  a non-orientable surface with punctures
is studied via classical homotopy theory of configuration spaces. In particular, we obtain a non-orientable version of the Birman exact sequence. In the case of $\rp^2$, we analize
the Serre spectral sequence of a fiber bundle $  \conf{2} \to X_k \to  BSO(3)  $
where $X_k$ is a $K(\Gamma^k(\rp^2),1)$ and  $\conf{2}$ denotes the configuration space of unordered $k$-tuples of distinct
points in $\rp^2$. As a consequence, we express the mod-2 cohomology of $\Gamma^k(\rp^2)$ in terms of that of
$\conf{2}$.
\end{abstract}



\section{Introduction}

Let $S_g$ be a closed orientable surface of genus $g$ and let $\diff^+(S_g)$ denote the group 
of orientation preserving diffeomorphisms of $S_g$, under the compact-open topology. 
The {\em mapping class group} $\Gamma_g$ is defined to be $\pi_0 \diff^+(S_g)$, the group of
path components for $\diff^+(S_g)$. Equivalently, $\Gamma_g$ can be defined as 
$\diff^+(S_g) / \diff_0(S_g)$, the group of isotopy classes of orientation preserving diffeomorphisms
of $S_g$, where $\diff_0(S_g)$ is the identity component of $\diff^+(S_g)$. This group plays an 
important role in the theory of Teichm\"uller spaces since $\Gamma_g$ acts on the space 
$\mathcal T_g$ of complex structures on $S_g$ and the quotient $\mathcal M_g$ is the moduli 
space of Riemann surfaces of genus $g$, see \cite{IT92}, \cite{FM11}. Moreover, its cohomology 
is closely related to the theory of characteristic classes of $S_g$-bundles, \cite{MO01}, \cite{EB08}.
There exist some variations on the definition above, for instance, consider the subgroup 
$\diff^+(S_g;k)$ of orientation preserving diffeomorphisms of $S_g$ which leave a fixed set of $k$ 
distinct points invariant. The mapping class group of $S_g$ with $k$ marked points $\Gamma_g^k$ 
is defined to be $\pi_0\diff^+(S_g;k)$, the group of path components for $\diff^+(S_g;k)$.
Similarly, if $N_g$ is a non-orientable surface of genus $g$, $\Gamma^k(N_g)$ is defined to be 
$\pi_0 \diff(N_g; k)$ where $\diff(N_g; k)$ is the group of diffeomorphisms of $N_g$ which leave 
a set of $k$ points invariant. In the special case when $k=0$ one recovers the 
{\em classical non-orientable} mapping class group $\Gamma(N_g) = \diff(N_g) / \diff_0(N_g)$.
\\

The purpose of this work is to use classical homotopy theory to compute the mod-2 cohomology of   
$\Gamma^k(\rp^2)$. Let $F_k(M)$ denote the configuration space of ordered $k$-tuples of distinct 
points in $M$, equipped with the natural action of the symmetric group $\Sigma_k$. Notice that 
$\diff(\rp^2)$ acts diagonally on the space $F_k(\rp^2)$ and thus it acts on the unordered configuration 
space $F_k(\rp^2)/\Sigma_k$. So one is led to consider the Borel construction

$$  E \diff(\rp^2)  \underset{\diff(\rp^2)}\times  \conf{2}.$$

\n
Using the classical result that  $\diff(\rp^2)$ has $SO(3)$ as a deformation retract \cite{AG73}, 
\cite{EE69}, the construction above is homotopy equivalent to

$$  E SO(3)  \underset{SO(3)}\times  \conf{2}.$$

\n
Moreover, it was proven in \cite{W02} that this space is an Eilenberg-MacLane space $K(\pi,1)$
for $k \geq 2$, and we show here that its fundamental group is isomorphic to $\Gamma^k(\rp^2)$. 
Therefore, one may use the universal fibration

\begin{equation}\label{E:mainfibration}
 \conf{2}  \longrightarrow  E SO(3)  \underset{SO(3)}\times  \conf{2}    \longrightarrow   B SO(3)
\end{equation}

\n
to study the cohomology of $\Gamma^k(\rp^2)$. In fact we show in Section 5 that the 
associated Serre spectral sequence in mod-2 cohomology collapses at the $E_2$-term. 
As a consequence, there is an isomorphim of   $H^*(BSO(3); \F_2)$--modules

$$   H^*(\Gamma^k(\rp^2); \F_2)  \cong  H^*(BSO(3); \F_2)  \otimes  H^*(\conf{2}; \F_2).$$

\n 
The additive structure of the mod-2 cohomology of $F_k(M)/\Sigma_k$ is well understood 
for $M$ a compact smooth manifold and it is determined by the dimension of $M$ and its 
Betti numbers, see \cite{BO89}. In the case of the projective plane, one can also express 
the mod-2 (co)homology of $\conf{2}$ in terms of that of the classical braid groups  $B_k$, 
see section 6.\\

Mapping class groups have also been exhaustively studied due to its relations to low 
dimensional topology. One of these relations is the classical result that the mapping class 
group of the 2-disk with $k$ marked points is isomorphic to $B_k$, the Artin's braid group on 
$k$-strands (\cite{FM11}). In the case of the sphere $S^2$, the group $\Gamma^k_0$ is 
isomorphic to the quotient of the braid group $B_k(S^2)$ by its center (\cite{BJ69}). In section 
3 we restrict attention to non-orientable surfaces $N_g$ and consider the $\diff_0(N_g)$-Borel 
construction  on the unordered configuration space
$$E \diff_0(N_g)  \underset{\diff_0(N_g)}{\times}   F_k(N_g)/\Sigma_k .  $$

\n
We define the {\em reduced mapping class group}  $\widetilde\Gamma^k(N_g)$ to be the 
fundamental group of this space. We show that $\widetilde\Gamma^k(N_g)$ is naturally a 
subgroup of $\Gamma^k(N_g)$ and fits into the following exact sequence

$$ 1  \longrightarrow  \widetilde\Gamma^k(N_g)  \longrightarrow  \Gamma^k(N_g)  
\longrightarrow  \Gamma(N_g)  \longrightarrow 1 .$$

\n
It is easy to show that for $g \geq 3$, the group $\widetilde\Gamma^k(N_g)$ is isomorphic 
to the surfce braid group $B_k(N_g)$. Thus, the exact sequence above recovers the 
{\em Birman exact sequence} in the non-orientable case. In the case of the real proyective 
plane $\rp^2$, the center of $B_k(\rp^2)$ is $\Z_2$ and we show there is an isomorphism

$$ \Gamma^k(\rp^2)  \cong  B_k(\rp^2) / \Z_2 . $$

\n
Similarly, in the case of the Klein bottle $\K$, we show the group $\widetilde\Gamma^k(\K)$ 
is isomorphic to $B_k(\K)$ modulo a central $\Z$.\\

The article is organized as follows. In section 2 we show how to construct concrete
$K(\pi,1)$ spaces for mapping class groups using configuration spaces. In section 3 
we study the group $\widetilde\Gamma^k(M) = \pi_0( \diff_0(M)  \cap  \diff(M;k) )$  
and derive the Birman exact sequence in the non-orientable case. In section 4 we 
show the Serre spectral sequence of bundle (\ref{E:mainfibration}) collapses and prove 
the main theorem. In section 5 we recall the main results from \cite{BO89} on the mod-2 
homology of configuration spaces for compact smooth manifolds, and finally, in section 6 
we carry out some explicit calculations for the homology of $F_k(M)/\Sigma_k$ in the 
case when $M=\rp^2$.\\


 \section{Configuration spaces and mapping class\\ groups}

Given a manifold $M$  and $k \geq 1$, define the configuration space of ordered $k$-tuples 
of distinct points in $M$ by

$$ F_k(M)  =\{(m_1, \dots, m_k)   \in M^k   \mid  m_i\neq m_j \;\; \text{if}  \;\; i\neq j  \} .$$

\n
The symmetric group on $k$ letters, $\Sigma_k$, acts naturally on $F_k(M)$ by permutation 
of coordinates and the orbit space $F_k(M)/\Sigma_k$ is the {\em unordered} configuration 
space. It is a classical result that when $M$ is a surface $\neq S^2, \rp^2$ the spaces 
$F_k(M)/\Sigma_k$ are Eilenberg-MacLane spaces $K(\pi,1)$ for the corresponding surface braid 
groups on $k$ strands, $B_k(M) = \pi_1 F_k(M)/\Sigma_k$. In the exceptional cases when 
$M=S^2$ or $\rp^2$, one needs to consider the Borel construction with respect to the natural 
$S^3$ action. Namely, let $G$ be a topological group, $EG$ a contractible space with a right free 
$G$-action and $X$ a left $G$-space. The associated homotopy orbit space (the Borel construction) 
$EG\times_G X$ is defined as the quotient of $EG\times X$ under the action of $G$ given 
by $g(x, m) = (xg^{-1}, gm)$. The projection onto the first coordinate induces a fiber bundle

$$ X \longrightarrow EG \times_G X   \longrightarrow  BG  $$

\n
where $BG$ denotes the classifying space of $G$.\\

Now, in the case when  $M = S^2$ or $\rp^2$, the corresponding $K(\pi,1)$'s for the braid 
groups $B_k(S^2)$ and $B_k(\rp^2)$ are given as the Borel constructions 
$ES^3 \times_{S^3} F_k(S^2)/\Sigma_k$ and $ES^3 \times_{S^3} F_k(\rp^2)/\Sigma_k$, 
respectively, where $S^3$ acts on the configuration spaces $F_k(S^2)/\Sigma_k$ and 
$F_k(\rp^2)/\Sigma_k$ via the double cover $S^3\to SO(3)$, see \cite{CO93}, \cite{W02}.
Some further Borel constructions provide $K(\pi,1)$ spaces for groups related to the mapping 
class groups as shown next. A proof of the following result can be found in \cite{tD08}, 
Theorem 1.8.6.\\

\begin{lemma}\label{hisot}
Let $G$ be a locally compact Hausdorff topological group with a countable basis acting 
transitively on a locally compact Hausdorff space $X$. Then for each $x\in X$ the map 
$G \to X$, given by $g \mapsto gx$ is open and the induced map
$$ G/G_x  \to X  $$
is a homeomorphism, where $G_x$ is the isotropy group of $x$.
\end{lemma}

\bigskip

\begin{corollary}
Under the assumptions of Lemma \ref{hisot} above, for every $x\in X$ there is a 
homotopy equivalence  $EG \times_G X  \simeq  B G_x$.
\end{corollary}

\medskip

\n
{\em Proof:} Notice that 
$EG \times_G X \approx EG \times_G G/G_x = EG/G_x \simeq  BG_x$.
$\square$ \\

\bigskip

\n
Now, let $M$ be a closed orientable surface and notice that $\diff^+(M)$ acts diagonally on 
$F_k(M)$ and thus on the unordered configuration space $F_k(M)/\Sigma_k$. The isotropy 
group of a fixed configuration  $Q_k = \{m_1, \dots, m_k\}$ in  $F_k(M)/\Sigma_k$ is 
$\diff^+(M;k)$, the group of orientation preserving diffeomorphisms of $M$ which leave the 
set $Q_k$ invariant. Thus, by the corollary above, there is a homotopy equivalence

$$ E \diff^+(M)   \underset{\diff^+(M)}\times F_k(M)/\Sigma_k \;  \simeq  \; B\diff^+(M;k) $$

\n
which induces an isomorphism of fundamental groups

\begin{align*}
\pi_1 \left( E \diff^+(M)   \underset{\diff^+(M)}\times F_k(M)/\Sigma_k \right)
\;  & \cong  \;  \pi_1 B\diff^+(M;k)  \\
&  \cong  \; \pi_0 \diff^+(M;k) .
\end{align*}

\n
Notice the last group is $\Gamma^k(M)$, the mapping class group of $M$ with $k$ 
marked points.\\

\bigskip

\n
{\bf Example:}
If $M = S^2$, a classical result of Smale \cite{Sm59} states that the inclusion 
$SO(3) \to \diff^+(S^2)$ is a homotopy equivalence. Thus the natural map

$$  
ESO(3)   \underset{SO(3)}\times  F_k(S^2)/\Sigma_k
\; \longrightarrow \;
E\diff^+(S^2)   \underset{\diff^+(S^2)}\times   F_k(S^2)/\Sigma_k 
$$
is  a homotopy equivalence and the fundamental group of this space is the 
mapping class group $\Gamma^k(S^2)$. Moreover, F.~Cohen proved in \cite{CO93}
 that for $k \geq 3$, the above construction is an Eilenberg-MacLane space   
$K(\Gamma^k(S^2),1)$  and the cohomology $H^*(\mcg^k(S^2); \mathbb F_2)$ 
with mod-2 coefficients was described by C.F.~B\"odigheimer, F.~Cohen and 
D.~Peim \cite{BCP01}.\\


 \section{The reduced mapping class group}


In this section $M$ will denote a non-orientable surface, although the same arguments 
apply to the orientable case with the obvious modifications. In the non-orientable case, 
the mapping class groups $\Gamma(M)$ and $\Gamma^k(M)$ are defined using the 
group $\diff(M)$  of all diffeomorphisms of $M$. The homotopy type of $\diff_0(M)$, the 
group of diffeomorphisms isotopic to the identity, is also known in this case, 
see \cite{EE69} and \cite{AG73}. We recall the result here:\\

\begin{theorem}[\cite{EE69}, \cite{AG73}]\label{eegram}
Let $N$ be a closed non-orientable surface, and let $\diff_0(N)$ be the group of 
diffeomorphisms isotopic to the identity. Then,
\begin{enumerate}
\item If $N=\rp^2$, then $\diff_0(N) =\diff(N)$ is homotopy equivalent to $SO(3)$.
\item If $N$ is the Klein bottle, then $\diff_0(N)$ is homotopy equiva\-lent to $SO(2)$.
\item If $N=N_g$,  a closed non-orientable surface of genus $g\geq 3$, $\diff_0(N)$ is contractible.
\end{enumerate}
\end{theorem}

\bigskip

Next we consider the $\diff_0(M)$-Borel construction. Notice $\diff_0(M)$ acts on the 
configuration space  $F_k(M)/\Sigma_k$ with isotropy subgroup  $\diff_0(M) \cap \diff(M;k)$.
Therefore, there is a homotopy equivalence

$$
E\diff_0(M)\underset{\diff_0(M)}{\times}F_k(M)/\Sigma_k\;\simeq \;B\left(\diff_0(M)\cap \diff(M;k)\right)
$$

\n
and thus and isomorphism of groups

\begin{align*}
\pi_1\left(E\diff_0(M)\underset{\diff_0(M)}{\times}F_k(M)/\Sigma_k\right)   &  
                       \cong \pi_1B \left(\diff_0(M)\cap \diff(M;k)\right)\\
& \cong  \pi_0\left(\diff_0(M)\cap \diff(M;k)\right)  .
\end{align*}

\n
Inspired by this isomorphism, we define the {\em reduced mapping class group} of $M$ with $k$ 
marked points by

$$
\widetilde{\Gamma}^k(M)=\pi_0( \diff_0(M) \cap     \diff(M;k))
$$

\bigskip

\n
This group is closely related to the {\it extended} and punctured mapping class group, 
as shown next.\\

\begin{theorem}\label{rmcgs}
Let $M $ be a compact connected surface. The reduced mapping class group 
$\widetilde{\mcg}^k(M)$ is naturally a subgroup of $\mcg^k(M)$ and fits into the 
following exact sequence:

$$
1\longrightarrow \widetilde{\mcg}^k(M) \longrightarrow \mcg^k(M)\longrightarrow \mcg(M)\longrightarrow 1.
$$
\end{theorem}

\bigskip
 
To prove this result, we will need a couple of simple lemmas.\\


\begin{lemma}
For every $k$, 
$\diff(M;k)  \cdot \diff_0(M) =  \diff(M)  $.
\end{lemma}

\n
{\em Proof:} Let $\varphi \in \diff(M)$ and let $\{x_1, \dots, x_k\}$ be a fixed configuration 
of $k$ distinct points in $M$. It is easy to see there is diffeomorphism  $\psi: M \to M$, 
isotopic to the identity, sending the points  
$\varphi^{-1}(x_1), \dots, \varphi^{-1}(x_k)$ to $x_1, \dots, x_k$. Then   
$  \varphi  = (  \varphi \circ \psi^{-1} ) \circ \psi$, where  $ \varphi \circ \psi^{-1} \in \diff(M;k)$  
and $\psi \in \diff_0(M)$. $\square$\\


\bigskip

\begin{lemma}\label{ll1}
For $k \geq 0$, the quotient group
$\diff(M;k)/\diff(M;k)\cap \diff_0(M)$ is isomorphic to $\Gamma(M)$.
\end{lemma}

\n
{\em Proof: } Consider the natural projection $p:\diff(M)\to \diff(M)/\diff_0(M)$ and let
$p_k$ be the restriction to the subgroup $\diff(M;k)$. Notice that the image of $\diff(M;k)$ 
is  isomorphic to the quotient

$$
\frac{\diff(M;k)\cdot\diff_0(M)}{\diff_0(M)}
$$

\n
and the kernel of $p_k$ is  $\diff(M;k)\cap \diff_0(M)$. Thus $p_k$ induces an isomorphism 
$$
\frac{\diff(M;k)}{\diff(M;k)\cap \diff_0(M)}   \xrightarrow{\;\cong\;} \frac{\diff(M;k)\cdot\diff_0(M)}{\diff_0(M)}
= \frac{\diff(M)}{\diff_0(M)}.
$$

$
\qquad
\qquad
\qquad
\qquad
\qquad
\qquad
\qquad
\qquad
\qquad
\qquad
\qquad
\qquad
\qquad
\qquad
\qquad
\square
$\\

\n
{\em Proof of Theorem \ref{rmcgs}:}
Let $\diff_0(M;k)$ be the connected component of the identity of $\diff(M;k)$ and 
consider the following diagram of short exact sequences of topological groups:

$$
\xymatrix{\diff_0(M;k)\ar@{=}[r]\ar[d]&\diff_0(M;k)\ar[d]\ar[r]&\{id\} \ar[d]\\\diff(M;k)\cap \diff_0(M)\ar[d]\ar[r]&\diff(M;k)\ar[d]\ar[r]&\frac{\diff(M;k)}{\diff(M;k)\cap \diff_0(M)}\ar[d]\\\frac{\diff (M;k)\cap \diff_0(M)}{\diff_{0} (M;k)}\ar[r]&\frac{\diff(M;k)}{\diff_0(M;k)}\ar[r]&\frac{\diff(M;k)}{\diff(M;k)\cap \diff_0(M)}}
$$

\medskip

\n
where the bottom row is an exact sequence of discrete groups. By Lemma \ref{ll1} the 
cokernel is isomorphic to $\Gamma(M)$. On the other hand, notice the path component 
of the identity of the group  $\diff_0(M) \cap \diff(M;k)$ is precisely $\diff_0(M;k)$. Then the 
kernel is given by

$$
\frac{\diff(M;k)\cap \diff_0(M)}{\diff_0(M;k)} \cong   \pi_0(\diff(M;k)\cap \diff_0(M))  = \tilde{\Gamma}^k(M)
$$
and the theorem follows. $\square$\\


We proceed to give some examples in low genus.  \\

\bigskip

\n
{\bf Example:} 
Let $M= \rp^2$ and consider the natural action of $SO(3)$ on $\rp^2$  given by rotation 
of lines through the origin in $\R^3$. It was shown by Gramain \cite{AG73} that the 
natural inclusion $SO(3) \to  \diff(\rp^2)$ is a homotopy equivalence, and thus
$\diff_0(\rp^2) = \diff(\rp^2)$ and $\widetilde\Gamma^k(\rp^2)  = \Gamma^k(\rp^2)$. 
Moreover, the natural map

$$ E SO(3)\times_{SO(3)} F_k(\rp^2)/\Sigma_k   \;  \xrightarrow{\;\; \simeq \;\;}  \;
E\diff(\rp^2) \times_{\diff(\rp^2)} F_k(\rp^2)/\Sigma_k  $$

\n
is a homotopy equivalence and the fundamental group of this space is $\Gamma^k(\rp^2)$.
It was shown in \cite{W02} that for $k \geq 2$ the Borel construction above is a $K(\pi,1)$. 
Thus:

\begin{theorem}
If $k\geq 2$, the $SO(3)$-Borel construction
$$ ESO(3) \times_{SO(3)}  F_k(\rp^2)/\Sigma_k $$
is an Eilenberg-MacLane space $K(\pi,1)$ where $\pi=\Gamma^k(\rp^2)$.
\end{theorem}

\bigskip

\n
On the other hand, notice there is a fibration
$$  B\Z_2 \longrightarrow  ES^3\times_{S^3} F_k(\rp^2)/\Sigma_k  \longrightarrow
ESO(3)\times_{SO(3)} F_k(\rp^2)/\Sigma_k $$

\n
of $K(\pi,1)$ spaces and thus one gets a short exact sequence of groups

$$ 1 \longrightarrow  \Z_2  \longrightarrow  B_k(\rp^2)  \longrightarrow  \Gamma^k(\rp^2)  \longrightarrow 1   $$

\medskip

\n
Moreover, recall that the center of $B_k(\rp^2)$ can be described as the image 
of the full twist braid $\Delta$ in $B_k(D^2)$ under the homomorphism   
$B_k(D^2)  \to  B_k(\rp^2)$ relative to an embedded disc $D^2 \subset  \rp^2$. 
This braid is known to be the only element of order 2, see \cite{GG10}. \\

\begin{corollary}
There is an isomorphism    $\Gamma^k(\rp^2)  \cong B_k(\rp^2)/\Z_2$,
where $\Z_2$ is the center of $B_k(\rp^2)$.
\end{corollary}

\bigskip

\n
{\bf Remark:}  Consider the matrices $\alpha, \beta, \gamma \in SO(3)$ given by
$$
\alpha=
\begin{pmatrix}
-1  &  0  &  0\\
0  &  1  &  0\\
0  &  0  & -1
\end{pmatrix}, \; \quad
\beta=
\begin{pmatrix}
1  &  0  &  0\\
0  &  -1  &  0\\
0  &  0  & -1
\end{pmatrix}, \; \quad
\gamma=
\begin{pmatrix}
1  &  0  &  0\\
0  &  0  &  -1\\
0  &  1  & 0
\end{pmatrix}.
$$
It can be shown that the group generated by $\alpha, \beta$ and $\gamma$ is 
isomorphic to $D_8$, the dihedral group of order 8. Moreover, it is proven in 
\cite{W02} that $SO(3)/D_8$ is homotopy equivalent to the unordered configuration 
space  $F_2(\rp^2)/\Sigma_2$. Thus, there is a homotopy equivalence

$$   E SO(3)   \underset{SO(3)}\times  SO(3)/D_8   \;  \simeq  \;
 E SO(3)   \underset{SO(3)}\times  F_2(\rp^2)/\Sigma_2  .
 $$

 \n
Since the space on the left is homotopy equivalent to the classifying space $BD_8$, 
we have that  $E SO(3)   \underset{SO(3)}\times  F_2(\rp^2)/\Sigma_2 $  is a $K(D_8,1)$
space.\\


\bigskip

\n
{\bf Example:} If $M= \K$ is the Klein bottle, it is well known that 
$\Gamma(\K) = \Z_2 \times \Z_2$, see \cite{Ham65}. Then the  exact sequence 
in Theorem \ref{rmcgs} is given by  
$$ 1 \longrightarrow \widetilde\Gamma^k(\K)  \longrightarrow  \Gamma^k(\K) 
\longrightarrow   \Z_2 \times \Z_2 \longrightarrow1$$

\n
and this exhibits $\widetilde\Gamma^k(\K)$ as a normal subgroup of index 4 of 
$\Gamma^k(\K)$. In this case  $SO(2) \simeq \diff_0(\K)$ and the natural map \\

$$ E SO(2)\underset{SO(2)}{\times} F_k(\K)/\Sigma_k   \;  \xrightarrow{\; \simeq \;}  \;
E\diff_0(\K) \underset{\diff_0(\K)}{\times} F_k(\K)/\Sigma_k  $$

\medskip

\n
is a homotopy equivalence. Also, it is easy to show this space is a $K(\pi,1)$. \\

\bigskip

\begin{theorem}
If $k\geq 1$, the $SO(2)$-Borel construction
$$ ESO(2) \times_{SO(2)}  F_k(\K)/\Sigma_k $$
is an Eilenberg-MacLane space $K(\pi,1)$ where $\pi=\widetilde\Gamma^k(\K)$.
\end{theorem}

\bigskip
\bigskip

\n
Finally, consider the natural fibration 

$$ F_k(\K)/\Sigma_k   \longrightarrow    ESO(2) \times_{SO(2)}  F_k(\K)/\Sigma_k  \longrightarrow  BSO(2) $$

\bigskip

\n
whose homotopy exact sequence is given by

$$  
1  \longrightarrow  \Z  \longrightarrow  B_k(\K)  \longrightarrow  \widetilde\Gamma^k(\K) \longrightarrow 1.
$$

\n
Then we have

\begin{corollary}
There is an isomorphism   $\widetilde\Gamma^k(\K)  \cong  B_k(\K)/\Z$.
\end{corollary}

\bigskip

\n
{\bf Remark:}
The subgroup $\Z$ also corresponds to the center of the braid group of $\K$, 
a fact that can be checked by two different methods. The first one considers an 
inductive argument using the Fadell-Neuwirth fibrations to relate the center of the 
pure braid group to the center $B_2(\K)$; then one uses the fact that any free group 
on at least two generators has trivial center. The other method consists of a direct 
calculation of the center from a given finite presentation of $B_k(\K)$. We thank 
F. Cohen and D. Gon\c{c}alves for nice talks dicussing these methods.
\\

\bigskip

\n
{\bf Example:} If $M=N_3$ it also well known that $\Gamma(N_3) = SL(2,\Z)$, 
see \cite{GB06}. Therefore, we have an exact sequence

$$ 1 \longrightarrow \widetilde\Gamma^k(N_3)  \longrightarrow  \Gamma^k(N_3) 
\longrightarrow  SL(2,\Z)\longrightarrow1$$
which shows that $\widetilde\Gamma^k(N_3)$ is a much smaller group than 
$\Gamma^k(N_3)$. In fact, we can prove the following result.\\

\begin{theorem}
For $g \geq 3$, the reduced mapping class group $\widetilde\Gamma^k(N_g)$ 
is isomorphic to the braid group  $B_k(N_g)$.
\end{theorem}

\medskip

\n
{\em Proof:}
Recall the group 
$\widetilde\Gamma^k(N_g) =  \pi_0( \diff_0(N_g) \cap  \diff(N_g; k)) $
is the fundamental group of the Borel construction

$$E \diff_0(N_g)  \underset{\diff_0(N_g)}{\times}   F_k(N_g)/\Sigma_g   $$

\n
On the other hand, projection onto the first coordinate induces a universal bundle
of the form

$$ 
F_k(N_g)/\Sigma_g \longrightarrow
E \diff_0(N_g)  \underset{\diff_0(N_g)}{\times}   F_k(N_g)/\Sigma_g   \longrightarrow
B \diff_0(N_g)$$ 

\n
But for $g \geq 3$, Theorem \ref{eegram} implies the classifying space $ B \diff_0(N_g) $ 
is contractible and the result follows. $\square$\\

\bigskip

Thus, Theorem \ref{rmcgs} recovers a version of the {\em Birman exact sequence} 
on the non-orientable case, see \cite{Bir69}, \cite{FM11}.\\


\bigskip

 \section{The mod-2 cohomology of $\Gamma^k(\rp^2)$}

\n
Recall that $SO(3)$ acts  on $\rp^2$ by rotating lines in $\R^3$:
$$ \mu:   SO(3)   \times   \rp^2  \longrightarrow   \rp^2$$
$$ \hspace{1in} ( A,    [x,y,z]     ) \quad \longmapsto \quad  [A(x,y,z)]$$
and thus, it acts diagonally on $\conf{2}$. So one may consider the 
$SO(3)$-Borel construction and the associated fibration

$$\conf{2}  \longrightarrow  ESO(3) \underset{\so}{\times} \conf{2} \longrightarrow  BSO(3).  $$

\bigskip

\noindent
\begin{theorem}\label{T:maintheorem}
The Serre spectral sequence for the fibration above collapses at the $E_2$-term,  
in mod 2 cohomology.
\end{theorem}

\medskip

\noindent
{\em Proof:}  
Notice the $SO(3)$ action on $\rp^2$ can be extended to an action on $\rp^\infty$  by setting:
$ \; \mu(A,  [x_0, x_1, x_2, x_3, \dots ]) =   [A(x_0, x_1, x_2), x_3, \dots]$. Thus the inclusion
$\conf{2}  \hookrightarrow \conf{\infty}$  is $SO(3)$-equivariant and gives rise to a map of 
fibrations

$$
\xymatrix
@C=2pc@R=4pc
{
\conf{2}  \ar[d] \ar@{^{(}->}[r]  &   \conf{\infty} \ar[d]    \\
\qquad ESO(3)    \underset{\so}{\times}  \conf{2}   \ar[d]    
& \qquad  {ESO(3)}   \underset{\so}{\times}  \conf{\infty}   \ar[d]     \\
  BSO(3)   \ar@{=}[r]  &   BSO(3)   \ar@<1ex> [ul]; [u]-<2cm,0cm>    
  }
$$

\n  
which induces a map between the corresponding spectral sequences. We will show in 
Theorem \ref{S:Handles} that the induced map on the fibers is an epimorphism in mod-2 
cohomology. Then the desired spectral sequence collapses provided the spectral sequence 
for the fibration on the right collapses. Secondly, notice the natural maps

$$     \conf{\infty}      \longleftarrow    E\Sigma_k  \underset{\;\Sigma_k}{\times}  \Conf{\infty}    
\longrightarrow E\Sigma_k  \underset{\;\Sigma_k}{\times} (\rp^\infty)^k$$
are $SO(3)$-equivariant  and also homotopy equivalences. Therefore we get equivalences of fibrations \\

\bigskip

{\scriptsize

\hspace{-1in}
$
\xymatrix@R=4pc
@C=0pc
{ 
\conf{\infty} \ar[d] &   E\Sigma_k \underset{\;\Sigma_k}{\times}  \Conf{\infty} \ar[d]  \ar@<-1ex>[l]_{\simeq} \ar@<1ex>[r]^{\simeq}&  
E\Sigma_k \underset{\;\Sigma_k}{\times}  (\rp^\infty)^k \ar[d]\\
\qquad  {ESO(3)}   \underset{\so}{\times}  \conf{\infty}   \ar[d]    
 &  
\quad \qquad \qquad ESO(3) \underset{\so}{\times}  \left[ E\Sigma_k \underset{\;\Sigma_k}{\times}  \Conf{\infty} \right]  \ar[d]
 &  
 \quad \qquad \qquad ESO(3)  \underset{\so}{\times}  \left[ E\Sigma_k \underset{\;\Sigma_k}{\times}  (\rp^\infty)^k \right]   \ar[d]\\
  BSO(3)   \ar@{=}[r]  &   BSO(3)  \ar@{=}[r]  &  BSO(3)    
 \ar[ul]; [u]-<1.7cm,0cm>
 \ar[ul]-<2cm,0,cm>; [ull]
}
$
}

\bigskip

\n
To prove the assertion of the theorem we will actually show the right column is a 
trivial fibration. Let $(\rp^\infty)_t$ denote the space $\rp^\infty$ endowed with the 
trivial $SO(3)$-action and consider the map 
$s:(\rp^\infty)_t  \to   \rp^\infty$
given by shifting of coordinates:  $s[x_0, x_1, x_2, \dots ]  = [0, 0, 0, x_0, x_1, x_2, \dots]$. 
Notice the map $s$ is $SO(3)$-equivariant and also a homotopy equivalence, since 
it is non-trivial on fundamental groups. Thus we get an equivalence of fibrations

\bigskip

{\small
\hspace{-1in}
$
\xymatrix@C=2pc@R=4pc{
 E\Sigma_k \underset{\;\Sigma_k}{\times}  (\rp^\infty)_t^k   \ar[d]\ar[r]^{1 \times s^k}_{\simeq}  &  
E\Sigma_k \underset{\;\Sigma_k}{\times}  (\rp^\infty)^k \ar[d]\\
  \quad \qquad \qquad ESO(3)  \underset{\so}{\times}  \left[ E\Sigma_k \underset{\;\Sigma_k}{\times}  (\rp^\infty)_t^k \right] \ar[d] 
 &  
 \quad \qquad \qquad ESO(3)  \underset{\so}{\times}  \left[ E\Sigma_k \underset{\;\Sigma_k}{\times}  (\rp^\infty)^k \right]   \ar[d]\\
 BSO(3)  \ar@{=}[r]  &  BSO(3)   
 \ar[ul];[u]-<1.9cm,0cm> ^(.75){\simeq}
 }
 $
 }

\bigskip

\n
But the left column of the previous diagram is clearly a trivial fibration, since the 
$SO(3)$-action on the fiber was trivial, and the statement of the theorem follows. 
$\square$\\

\bigskip

As a consequence of Theorem \ref{T:maintheorem} we have

\begin{theorem}
For $k \geq 2$, there is an isomorphism of $H^*(BSO(3);\F_2)$-modules

$$ H^*(ESO(3) \underset{\so}{\times} \conf{2} ; \F_2)  \cong  \F_2[w_2, w_3] \otimes H^*(\conf{2} ; \F_2)$$
where $w_2, w_3$ are the Stiefel-Whitney classes in the cohomology of $BSO(3)$.
\end{theorem}

\bigskip

\n
It follows from here that the mod-2 cohomology of $\Gamma^k(\rp^2)$ 
is determined by the mod-2 cohomology of $\conf{2}$.

\section{Homology of configuration spaces}

In this section we describe the homology of the space $\conf{2}$ in terms of the 
homology of a much larger space. Let $M$ be manifold of dimension $m$ and $X$ 
a connected space with base point $*$ and recall the labelled configuration space 
$C(M;X)$ is given by

$$ C(M;X) = \left(  \coprod_{k \geq 0}  F_k(M)  \underset{\Sigma_k}{\times}  X^k   \right)  \Big/  \approx$$

\n
where the relation $\approx$ is generated by
$$[m_1, \dots, m_k; x_1, \dots, x_k] \approx  [m_1, \dots, m_{k-1}; x_1, \dots, x_{k-1}] $$
if $x_k=*$. Such spaces of labelled configurations occur in \cite{BO85}, \cite{BO89}, 
\cite{MD75} and \cite{May72} as models for mapping spaces. The space $C(M;X)$ is 
naturally filtered by length of configurations 
$$  * = C_0(M;X)  \subseteq  C_1(M;X)  \subseteq \ldots  \subseteq   C(M;X)$$

\n
and the filtration quotients 
$C_k(M;X)  /  C_{k-1}(M;X)$  are denoted by $D_k(M;X)$.
The basic properties of this construction are given next, see \cite{May72}, \cite{BO89}, \cite{VS74}:

\begin{enumerate}
\item  If $M = \R^n$ then $C(\R^n;X)  \simeq \Omega^n\Sigma^n X$.

\item There is an analogue of Snaith's stable splitting
$$ C(M; X)   \;  \underset{s}{\simeq}  \;  \bigvee _{k=1}^\infty  D_k(M; X) .$$

\item There is a natural vector bundle over the configuration space
$$\eta_k :   F_k(M)  \times_{\Sigma_k} \R^k  \to  F_k(M)/\Sigma_k  $$
and the Thom space of its $n$-fold sum is homeomorphic to $D_k(M;S^n)$.
Therefore, by the Thom isomorphism 

$$ H_q  ( F_k(M)/\Sigma_k  )  \cong   \bar{H}_{q+kn}  D_k(M;S^n)     $$

\end{enumerate}

\begin{theorem}[\cite{BO89}]\label{T:BCT}
For $M$ a smooth, compact manifold of dimension $m$ and $X=S^n$, 
there is an isomorphism of graded vector spaces
$$ \theta  :  H_*(C(M; S^n) ;  \F_2) \; \cong \; 
\bigotimes_{q=0}^m  H_*(\Omega^{m-q}  S^{m+n} ; \F_2)^{\otimes\beta_q},  $$
where $\beta_q$ is the $q$-th Betti number of $M$.

\end{theorem}

\

Each factor   $H_*(\Omega^{m-q}S^{n+q})$  is an algebra with weights associated 
to its generators. This yields a filtration on the tensor product and  
$\bar{H}_* D_k(M; S^n)$ corresponds, via the isomorphism $\theta$, to the vector 
space generated by the elements of weight $k$. For completeness, we record here 
the mod 2 homology of the iterated loop spaces $\Omega^kS^{n+k}$, see \cite{FC87}. 
Throughout the rest of this and next section, all homology groups are taken with mod-2 
coefficients.  \\

Recall there are mod 2 homology operations which are natural for $n$-fold loop maps
$$  Q_i  :  H_q(\Omega^nX)  \longrightarrow  H_{2q+i} (\Omega^nX), \quad  0\leq i \leq n-1 $$
which are linear if $0 \leq i < n-1$, known as the {\em Dyer-Lashof operations}  
\cite{DL62}. Let $x_n \in H_n(\Omega^kS^{n+k})$ be the fundamental class and let 
$Q_I x_n$ denote the composition $Q_{i_1} Q_{i_2} \dots Q_{i_r} x_n$  if 
$I = (i_1, i_2, \dots, i_r)$. The sequence $I$ is admissible if $0 < i_1 \leq i_2 \leq \dots \leq i_r$; 
write $\lambda(I) \leq q$ if $i_r \leq q$ and $\ell(I)=r$. \\

\begin{theorem}[\cite{FC87}]\label{T:iterated}
There is and isomorphism of Hopf algebras
$$ H_*(\Omega^k S^{n+k})  \cong   \F_2[Q_Ix_n], \qquad n \geq 1, $$
for admissible $I$ with $\lambda(I) \leq k-1$  and $Q_I x_n$ is primitive.
\end{theorem}

\bigskip

Thus, to compute $H_*C(M; S^n)$, we first introduce for every basis element  
$\alpha \in H_q(M)$ a generator, namely the fundamental class
$u_\alpha   \in H_*(\Omega^{m-q} S^{m+n})$  of degree  $|u_\alpha| = q +n$ and weight  
$\omega(u_\alpha) =1$. Secondly, for each $u_\alpha$ and index  $I=(i_1,  i_2,\dots, i_r)$ 
there is an additional generator
$   Q_I  u_\alpha = Q_{i_1} Q_{i_2} \dots Q_{i_r} u_\alpha$ if the condition 
$0 < i_1 \leq i_2 \leq \dots \leq i_r < m-q  $ holds. We have:

\begin{itemize}
\item $|Q_I u_\alpha |  = i_1 + 2 i_2 + 4 i_3 + \dots  +  2^{r-1} i_r  + 2^r(q +n)$
\item  $\omega(Q_I u_\alpha)  =  2^{\ell(I)}  = 2^r$
\end{itemize}
and $u_\alpha^2=0$ if $|\alpha| = q = m$. \\

The isomorphism $\theta$ in Theorem \ref{T:BCT} depends on the choice of 
a handle decomposition for $M$ and it is natural for embeddings which respect 
the handle decompositions. Recall that a manifold $\bar{M}$ is obtained from a 
submanifold $M\subset \bar{M}$ of codimension $0$ by attaching a handle of 
index $q$ if $\bar{M} = M \cup D$ with $D \approx [0,1]^m$  and 
$M \cap D \approx [0,1]^{m-q} \times \partial [0,1]^q$. A handle decomposition 
for a manifold $M$ is a filtration by submanifolds

$$ 
M_0 \subset  M_1 \subset \ldots \subset  M_{m-1} \subset  M_m = M,
$$ 
where $M_q$ is obtained from $M_{q-1}$ by attaching handles of index $q$ and 
$M_0$ is a disjoint union of closed $m$-dimensional discs.\\

As an example, consider the usual CW structure on $S^m$ with two antipodal 
$q$-cells in every dimension, for $q = 0, 1, ..., n$, so that the $q$-skeleton of 
$S^m$ is the sphere $S^q$. It is clear that the cell structure above can be 
thickened to induce a handle decomposition for $S^m$ with {\it two antipodal handles} 
of index $q$, for $q=0, \dots, m$, such that the natural embedding 
$S^m \subset S^{m+1}$ respects the handle decomposition. \\

\begin{figure}[ht]
\vspace{-2.4in}
\includegraphics[scale=.6]{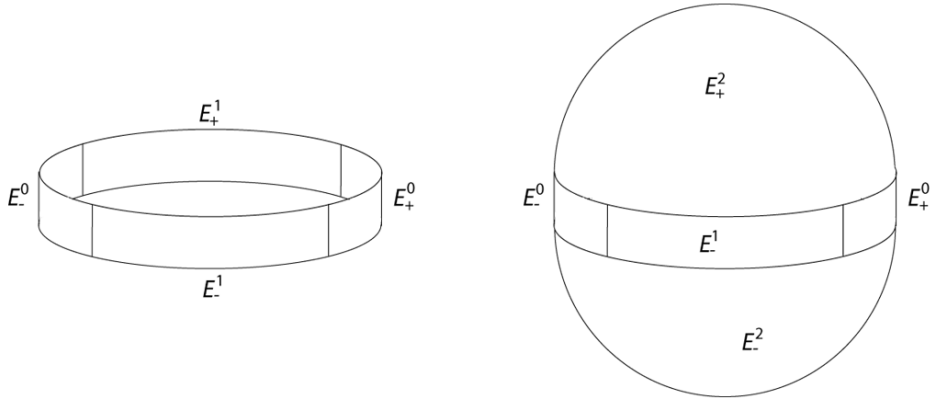}
\vspace{-2.4in}
\caption{Handle decompositions of $S^1$ and $S^2$}
\end{figure}

\n
Moreover, passing to the quotient by the antipodal $\Z_2$-action, we get a handle 
decomposition for $\rp^m$ with one handle of index $q$, for $q=0, \dots, m$, so 
that the natural embedding $\rp^m \subset \rp^{m+1}$ respects the corresponding 
handle decompositions. Thus, as a direct application of Theorem \ref{T:BCT}
we get \\

\begin{theorem}\label{S:Handles}
The natural inclusion $\rp^2 \subset \rp^\infty$ induces a map of unordered 
configuration spaces, which is an epimorphism in mod-2 cohomology:

$$ H^*(\conf{\infty}; \F_2)    \longrightarrow H^*(\conf{2}; \F_2) .$$  

\end{theorem}

\medskip

\n
{\em Proof:}
Let $m \geq 2$ and
consider the induced inclusion at the level of labelled configuration spaces  
$C(\rp^2; S^n)  \hookrightarrow  C(\rp^m; S^n)  $. By Theorem \ref{T:BCT} 
we have isomorphims in mod-2 cohomology:

\begin{align*}
H_* C(\rp^2; S^n)  & \cong  H_*(\Omega^2 S^{n+2})  \otimes H_*(\Omega S^{n+2}) \otimes H_*(S^{n+2}) \qquad \quad \text{and}\\
& \\
H_* C(\rp^m; S^n)   & \cong   H_*(\Omega^m S^{n+m})  \otimes H_*(\Omega^{m-1} S^{n+m})  \otimes 
\ldots \otimes H_*(S^{n+m}).
\end{align*}

\n
Moreover, the embedding  $\rp^2 \subset \rp^m$ respects the handle decompositions 
and the induced map in homology 
 $\; H_* C(\rp^2; S^n)  \to H_* C(\rp^m; S^n) \;$  
is given by the usual adjunction maps:
\begin{align*}
 & H_*(\Omega^2 S^{n+2})  \to   H_*(\Omega^m S^{n+m}) \\
 &  \\
 &  H_*(\Omega S^{n+2})  \to    H_*(\Omega^{m-1} S^{n+m}) \\
&  \\
& H_*(S^{n+2})  \to  H_*(\Omega^{m-2} S^{n+m}) 
 \end{align*}

\n
which are monomorphisms and preserve the weight of the generators. Therefore, the inclusion 
$\rp^2 \subset \rp^\infty$  induces a monomorphism in mod-2 homology

$$   H_*(F_k(\rp^2)/\Sigma_k; \F_2) \; \longrightarrow \;  H_*(F_k(\rp^\infty)/\Sigma_k;  \F_2)   $$

\n
and thus by duality, an epimorphism in mod-2 cohomology. $\square$  \\


\section{Explicit calculations}
Let us specialize to the case when $M$ is a surface. By Theorem \ref{T:BCT},  the mod-$2$ 
homology of the labelled configuration space $C(M;S^n)$ is given by the tensor product 

$$
H_*(\Omega^2S^{n+2})^{\otimes \beta_0}\otimes H_*(\Omega S^{n+2})^{\otimes \beta_1}
\otimes H_*(S^{n+2})^{\otimes \beta_2},
$$

\medskip

\n
where $\beta_0,\beta_1$ and $\beta_2$ are the mod-$2$ Betti numbers of $M$. In the case of 
the closed non-orientable surface $N_g$ of genus $g$ we have
$$  H_* C(N_g;S^n)  \; \cong \; \F_2[y_0,y_1,\ldots]\otimes \F_2[x_1,x_2,\ldots,x_g]\otimes \F_2[u]/u^2,  $$

\n
where $y_0, x_1, \ldots, x_g$ and $u$ are the fundamental classes on degrees $n,n+1$ and $n+2$, respectively, 
and $y_j=Q^j_1y_0=Q_1 Q_1 \ldots Q_1y_0$. Here $Q_1$ is the first Dyer-Lashof operation, 
so notice $|y_j|=(2^j-1)+2^jn$. The weights of all generators are given by
\begin{align*}
&\omega(u)=1,\\
&\omega(x_i)=1,\;\; \mbox{for} \;i=1,\ldots, g ,\\
&\omega(y_j)=2^j ,
\end{align*}

\n
and a basis for $\overline{H}_qD_k(N_g;S^n)$ consists of all monomials of the form
$$
h=u^ex_1^{a_1}\dots x_g^{a_{g}}y_0^{b_0}y_1^{b_1}\ldots y_r^{b_r},
$$
for some $r\geq 0,\;e=0,1$ and $a_i,b_j\geq 0$, such that
$$
\omega(h)=e+a_1+\ldots +a_g+b_0+2b_1+\ldots+2^rb_r=k
$$

\noindent
For example, $\overline{H}_qD_2(N_g;S^n)$ is determined by the following table.
\vspace{5mm}
\begin{table}[ht]
\centering  
\begin{tabular}{lll} 
\hline
$q$ & basis & rank \\ [0.5ex] 
\hline                  
$2n$ & $y_0^2$ & $1$   \\ 
$2n+1$ &$x_1y_0, \ldots , x_gy_0, y_1$ & $g+1$  \\
$2n+2$ & $uy_0, x_1^2, \ldots, x_ix_j, \ldots, x_g^2$ & $1+(g^2+1)/2$   \\
$2n+3$ & $ux_1, \ldots, ux_g$ & $g$ \\ [1ex]      
\hline 
\end{tabular}
\end{table}

\bigskip

\n
Thus, the rank of $H_q(F_2(N_g)/\Sigma_2;\F_2)$ is $1,g+1,1+(g^2+g)/2$ and $g$ for 
$q=0,1,2,3$. Similarly, in the case of  $\rp^2$, $g=1$, and  $k=3, 4$ one  has \\

$$ 
H_q(F_3(\rp^2)/\Sigma_3;\F_2)= 
\begin{cases}
\F_2 &\mbox{if}\;\;q=0\\
\F_2^2 &\mbox{if}\;\;q=1\\
\F_2^3 &\mbox{if}\;\;q=2 \\
\F_2^3 &\mbox{if}\;\;q=3 \\
\F_2 &\mbox{if}\;\;q=4\\
0&\mbox{otherwise}
\end{cases}
$$

\n
and

$$
H_q(F_4(\rp^2)/\Sigma_4;\F_2)=\begin{cases}
\F_2 &\mbox{if}\;\;q=0\\
\F_2^2 &\mbox{if}\;\;q=1\\
\F_2^4 &\mbox{if}\;\;q=2 \\
\F_2^5 &\mbox{if}\;\;q=3 \\
\F_2^3 &\mbox{if}\;\;q=4\\
\F_2 &\mbox{if}\;\;q=5\\
0&\mbox{otherwise}.
\end{cases} 
$$

\

\

\n

Finally, we express the mod-2 homology of $F_k(\rp^2)/\Sigma_k$ in terms of the 
homology of the classical braid groups. Since the Betti numbers of $\rp^2$ are 
$\beta_0=\beta_1=\beta_2 = 1$, then by Theorem 5.3 there is an isomorphism 

$$H_* C (\rp^2; S^n )  \; \cong  \; \F_2 [u]/u^2 \otimes \F_2[x] \otimes H_*(\Omega^2 S^{n+2} )  , $$

\bigskip

\n
where   $|u|= n + 2$,    $ |x| = n + 1$  and $\omega(u) = \omega(x) = 1$. 
Thus, a basis for $H_* C(\rp^2; S^n)$ is given by all the monomials of the form \\

 $$u^\epsilon \otimes x^\ell \otimes  y,   \qquad \quad  \epsilon =0,1,  \quad \ell =0, 1, 2, \dots $$

 \bigskip
 
 \n
 where $y$ runs over an additive basis for $H_*(\Omega^2S^{n+2})$. Notice that basis 
 elements of the form $u^\epsilon\otimes x^\ell$ have degree and weight given by

 \begin{align*}
  |u^\epsilon\otimes x^\ell|   & = \epsilon(n + 2) + \ell(n + 1) = \epsilon n + 2\epsilon + \ell n + \ell,\\
   & \\
  \omega(u^\epsilon\otimes x^\ell)      & =  \epsilon + \ell
 \end{align*}

\medskip
  
\noindent
Therefore  for fixed $q$ and  $k$, the monomial  $u^\epsilon\otimes x^\ell  \otimes y$  
represents a generator in 

$$ H_q  F_k(\rp^2)/\Sigma_k    \cong    H_{q + kn}   D_k(\rp^2; S^n)$$
  
\medskip
  
\n
if and only if $y$ has degree:
 \begin{align*}
|y|  &  = (q + kn) - (\epsilon n + 2\epsilon + \ell n+ \ell)  \\
      &  =  q + (k-\ell-\epsilon)n - 2\epsilon -\ell  
\end{align*}

\n
and weight:   
$ \omega(y) = k  - (\epsilon + \ell) =  k -\epsilon - \ell $.\\

\bigskip

\noindent
{\bf Case $\epsilon =0$:}
Elements in $ H_*(\Omega^2S^{n+2}) \cong  H_* C(\R^2; S^n)$ of degree  $q+ (k-\ell)n -\ell$ and weight  $k-\ell$
generate the vector space

\begin{align*} 
H_{q + (k-\ell)n -\ell}  \left( D_{k-\ell} (\R^2;  S^n)   \right)  &
\cong
H_{q-\ell} \left( F_{k-\ell}(\R^2)/\Sigma_{k-\ell}  \right)  \\
& =
H_{q-\ell}(B_{k-\ell})
\end{align*}

\bigskip

\n
{\bf Case $\epsilon =1$:}
Elements in $  H_* C(\R^2; S^n)$ of degree  $q+ (k-\ell-1)n -\ell-2$ and weight  $k-\ell -1$
generate the vector space

\begin{align*} 
H_{q + (k-\ell-1)n -\ell-2}  \left( D_{k-\ell-1} (\R^2;  S^n)  \right)   &
\cong
H_{q-\ell-2} \left( F_{k-\ell-1}(\R^2)/\Sigma_{k-\ell -1}  \right) \\
& =
H_{q-\ell-2}(B_{k-\ell-1})
\end{align*}


\medskip

\noindent
Thus there is an isomorphism of graded vector spaces:

$$  H_q  F_k(\rp^2)/\Sigma_k     \cong      \bigoplus_{\ell=0}^{\text{min}\{q,k\}} H_{q-\ell}(B_{k-\ell})   
\;  \oplus  \;    \bigoplus_{\ell=0}^{\text{min}\{q-2,k-1\}}  H_{q-\ell-2}(B_{k-\ell-1})$$

\bigskip

\n
It is worth to compare this result with the analog for $S^2$ which is obtained in \cite{BCP01} in a similar manner to the one exposed here.\\

\bibliographystyle{amsplain}
\bibliography{Paper2017}

\providecommand{\bysame}{\leavevmode\hbox to3em{\hrulefill}\thinspace}
\providecommand{\MR}{\relax\ifhmode\unskip\space\fi MR }
\providecommand{\MRhref}[2]{%
  \href{http://www.ams.org/mathscinet-getitem?mr=#1}{#2}
}
\providecommand{\href}[2]{#2}
\begin{thebibliography}{10}

\bibitem{Bir69}
J.S. Birman, \emph{Braids, links and mapping class groups}, Annals of
  Mathematical Studies, vol.~82, Princeton University Press, 1969.

\bibitem{BJ69}
\bysame, \emph{On braid groups}, Comm. Pure Appl. Math. \textbf{22} (1969),
  41--72.

\bibitem{BO85}
C.F. B\"odigheimer, \emph{Stable splittings of mapping spaces}, Algebraic
  topology ({S}eattle, {W}ash., 1985), Lecture Notes in Math., vol. 1286,
  Springer, Berlin, 1987, pp.~174--187.

\bibitem{BO89}
C.F. B{\"o}digheimer, F.~Cohen, and L.~Taylor, \emph{On the homology of
  configuration spaces}, Topology \textbf{28} (1989), no.~1, 111--123.

\bibitem{BCP01}
C.F. B\"odigheimer, F.R. Cohen, and M.D. Peim, \emph{Mapping class groups and
  function spaces}, Homotopy methods in algebraic topology ({B}oulder, {CO},
  1999), Contemp. Math., vol. 271, Amer. Math. Soc., 2001, pp.~17--39.

\bibitem{FC87}
F.R. Cohen, \emph{A course in some aspects of classical homotopy theory},
  Algebraic topology ({S}eattle, {W}ash., 1985), Lecture Notes in Math., vol.
  1286, Springer, Berlin, 1987, pp.~1--92.

\bibitem{CO93}
\bysame, \emph{On the mapping class groups for punctured spheres, the
  hyperelliptic mapping class groups, {${\rm SO}(3)$}, and {${\rm
  Spin}^c(3)$}}, Amer. J. Math. \textbf{115} (1993), no.~2, 389--434.
  \MR{1216436}

\bibitem{DL62}
E.~Dyer and R.K. Lashof, \emph{Homology of iterated loop spaces}, Amer. J.
  Math. \textbf{84} (1962), 35--88.

\bibitem{EE69}
C.~J. Earle and J.~Eells, \emph{A fibre bundle description of {T}eichm\"uller
  theory}, J. Differential Geometry \textbf{3} (1969), 19--43.

\bibitem{EB08}
J.~Ebert and O.~Randal-Williams, \emph{On the divisibility of characteristic
  classes of non-oriented surface bundles}, Topology Appl. \textbf{156} (2008),
  no.~2, 246--250.

\bibitem{FM11}
B.~Farb and D.~Margalit, \emph{A primer on mapping class groups}, Princeton
  Mathematical Series, vol.~49, Princeton University Press, 2012.

\bibitem{GG10}
D.L. Gon\c{c}alves and J.~Guaschi, \emph{Braid groups of non-orientable
  surfaces and the {F}adell-{N}euwirth short exact sequence}, J. Pure Appl.
  Algebra \textbf{214} (2010), no.~5, 667--677.

\bibitem{GB06}
F.J. Gonz{\'a}lez-Acu{\~n}a and J.M. M{\'a}rquez-Bobadilla, \emph{On the
  homeotopy group of the non orientable surface of genus three}, Rev.
  Colombiana Mat. \textbf{40} (2006), no.~2, 75--79.

\bibitem{AG73}
A.~Gramain, \emph{Le type d'homotopie du groupe des diff\'eomorphismes d'une
  surface compacte}, Ann. Sci. \'Ecole Norm. Sup. (4) \textbf{6} (1973),
  53--66.

\bibitem{Ham65}
M.E. Hamstrom, \emph{Homotopy properties of the space of homeomorphisms on
  {$P^{2}$} and the {K}lein bottle}, Trans. Amer. Math. Soc. \textbf{120}
  (1965), 37--45.

\bibitem{IT92}
Y.~Imayoshi and M.~Taniguchi, \emph{An introduction to {T}eichm\"uller spaces},
  Springer-Verlag, Tokyo, 1992.

\bibitem{May72}
J.P. May, \emph{The geometry of iterated loop spaces}, Lecture Notes in Math.,
  vol. 271, Springer-Verlag, Berlin-New York, 1972.

\bibitem{MD75}
D.~McDuff, \emph{Configuration spaces of positive and negative particles},
  Topology \textbf{14} (1975), 91--107.

\bibitem{MO01}
S.~Morita, \emph{Geometry of characteristic classes}, Translations of
  Mathematical Monographs, vol. 199, American Mathematical Society, 2001.

\bibitem{Sm59}
S.~Smale, \emph{Diffeomorphisms of the 2-sphere}, Proc. Amer. Math. Soc.
  \textbf{10} (1959), 621--626.

\bibitem{VS74}
V.P. Snaith, \emph{A stable decomposition of {$\Omega ^{n}S^{n}X$}}, J. London
  Math. Soc. (2) \textbf{7} (1974), 577--583.

\bibitem{tD08}
T.~tom Dieck, \emph{Algebraic topology}, EMS Texts in Mathematics, European
  Mathematical Society, 2008.

\bibitem{W02}
J.H. Wang, \emph{On the braid groups for {${\bf R}{\rm P}^2$}}, J. Pure Appl.
  Algebra \textbf{166} (2002), no.~1-2, 203--227. \MR{1868546}

\end{thebibliography}

\end{document}